\documentclass{article}%
\usepackage{amsfonts}
\usepackage{amsmath}
\usepackage{amsthm,amstext,amsfonts,bm,amssymb}
\usepackage{amssymb}
\usepackage{xfrac}
\usepackage{refcount}
\usepackage[margin=1cm,singlelinecheck=off]{caption}
\usepackage{enumitem}
\usepackage{graphicx}%
\usepackage{xypic}
\providecommand{\U}[1]{\protect\rule{.1in}{.1in}}
\providecommand{\comment}[1]{}

\newcommand{\CS}{\mathcal S}

\newtheorem{theorem}{Theorem}[section]

\newtheorem*{namedtheorem}{\theoremname}
\newcommand{\theoremname}{testing}

\makeatletter
\newtheorem*{rep@theorem}{\rep@title}
\newcommand{\newreptheorem}[2]{%
\newenvironment{rep#1}[1]{%
 \def\rep@title{#2 \ref{##1}}%
 \begin{rep@theorem}}%
 {\end{rep@theorem}}}
\makeatother

\newreptheorem{theorem}{Theorem}
\newreptheorem{lemma}{Lemma}
\newreptheorem{proposition}{Proposition}

\newtheorem{question}{Question}

\newtheorem{lemma}[theorem]{Lemma}

\newtheorem*{lemma*}{Lemma}

\newtheorem{remark}[theorem]{Remark}

\theoremstyle{definition}

\newtheorem*{definition*}{Definition}

\newcommand{\BZ}{\mathbb{Z}}
\newcommand{\BR}{\mathbb{R}}
\newcommand{\BQ}{\mathbb{Q}}
\newcommand{\BH}{\mathbb{H}}
\newcommand{\rank}{\mathrm{rank}}

\begin{document}

\title{Nielsen equivalence in mapping tori over the torus}
\author{Ian Biringer}
\maketitle

\begin{abstract}
We use the geometry of the Farey graph to give an alternative proof  of the fact that if $A \in GL_2\BZ$ and $G_A=\BZ^2 \rtimes_A \BZ$ is generated by two elements, there is a single Nielsen equivalence class of $2$-element generating sets for $G_A$ unless $A$  is conjugate to $\pm \left(\begin {smallmatrix} 2 & 1 \\ 1 & 1 \end {smallmatrix}\right )$, in which case there are two.

\vspace{2mm}

 \noindent \it MSC code: \rm 57M07. 
 
\noindent  \it Keywords: \rm Nielsen equivalence, Farey graph.
\end{abstract}
\section{Introduction}

\label {toruscase}

Let $G$ be a  finitely generated group. Two ordered $n $-element generating sets $S,T$ for $G $ are \emph {Nielsen equivalent} if the associated surjections $F_n \longrightarrow G$ differ by precomposition with a free group automorphism. This is equivalent to requiring that $S,T$ are related by a sequence of \emph{Nielsen moves}:\begin {enumerate}
\item if $a\neq b$ are generators, replace $a$ with $ab$,
\item if $a\neq b$ are generators, switch their places in the ordering.
\item if $a$ is a generator, replace it with $a^{-1}$,
\end {enumerate}
as the associated automorphisms generate $\mathrm{Aut}(F_n)$, see \cite[Chap. I, Prop. 4.1]{lyndon2015combinatorial}.

In \cite{Levittrank}, Levitt--Metaftsis studied Nielsen equivalence within groups of the form $G_A = \BZ^d \rtimes_A \BZ$,  where $A\in GL_d \BZ$.   Using the Cayley--Hamilton theorem, they show that $G_A$ is $2$-generated exactly when there is a vector $v\in \BZ^d$ such that $\langle v, Av \rangle = \BZ^d$.  They also show that the number of Nielsen equivalence classes of $2$-element generating sets is the index of $\langle A, -Id \rangle$ in its $GL_d \BZ$-centralizer. 

When $d=2$, one  can combine this with an observation of Cooper--Scharlemann \cite[Lemma 5.1]{cooper1999structure} to prove  the following theorem.

\begin {theorem}\label {torus2}
If $A \in GL_2\BZ $ and $G_A=\BZ^2 \rtimes_A \BZ$ is $2$-generated, there is a single Nielsen equivalence class of $2$-element generating sets for $G_A$ unless $A$  is conjugate to $\pm \left(\begin {smallmatrix} 2 & 1 \\ 1 & 1 \end {smallmatrix}\right )$, in which case there are two.
\end {theorem}
Note that when $A=\left(\begin {smallmatrix} 2 & 1 \\ 1 & 1 \end {smallmatrix}\right )$, $G_A$ is $2$-generated, since $\left<\left(\begin {smallmatrix} 1 \\ 0 \end {smallmatrix}\right ) ,\left(\begin {smallmatrix} 2 \\ 1 \end {smallmatrix}\right )\right > =\BZ^2$.  

\vspace{2mm}

 Our goal here is not to prove anything new, but rather to understand how to prove Theorem \ref{torus2} using the geometry of the \emph {Farey graph} $\mathcal F$. Algebraically, vertices of $\mathcal F$ are primitive elements $v=(p,q) \in \BZ^2$  up to negation, and vertices $v,w$ are connected by an edge if together they generate $\BZ^2$.  Any matrix $A \in GL_2\BZ $ acts on $\mathcal F$, and it turns out that Nielsen equivalence classes of $2$-element generating sets of $G_A$  correspond to geodesics in $ \mathcal F$ on which $A$ acts as a unit translation, see \S \ref{proof}.     Using this perspective, one can then prove Theorem \ref{torus2}  just using separation properties of geodesics in $ \mathcal F$.

 In the paper referenced above, Cooper--Scharlemann were interested in an analogue of  Theorem \ref{torus2} in the world of Heegaard splittings.  Recall that a closed surface $S$  in a  closed, orientable  $3$-manifold is a \emph{Heegaard splitting} if $M\setminus H$  has two components, each of which are (open) handlebodies. They showed  that  there is a unique  minimal genus Heegaard splitting of $M_A$  up to isotopy unless 
 $A$  is conjugate to $\pm \left(\begin {smallmatrix} 2 & 1 \\ 1 & 1 \end {smallmatrix}\right )$, in which case there are two.

Any  Heegaard splitting gives a pair of generating sets for $\pi_1 M$, just by taking free bases for the fundamental groups of the two handlebodies. These generating sets are well-defined up to Nielsen equivalence, and their Nielsen types certainly do not change if the Heegaard splitting $S$ is isotoped in $M$.  However,  in general it is hard to say when a generating set for $\pi_1 M$ is `geometric', i.e.\  when its Nielsen class comes from a Heegaard splitting, and when two (say, nonisotopic) Heegaard  splittings  give the same Nielsen class, see e.g.\ Johnson \cite{johnson2010horizontal}.

 However, inspired by the fact that the Cooper--Scharlemann result also applies when the minimal genus of a Heegaard splitting is $3$,  we ask:
 \begin {question}
 Is it true that if $\rank(G_A)=3$, there is a single Nielsen equivalence class of $3$-element generating sets?
\end {question}

Here, \emph{rank} is the minimal size of a generating set. In \cite{Biringerranks}, the author and Souto studied rank and Nielsen equivalence for mapping tori $M_\phi$, where $\phi : S \longrightarrow S$  is a pseudo-Anosov homeomorphism of a  closed orientable surface of genus $g\geq 2$. We showed that \emph {as long as $\phi$ has large translation distance in the curve complex $C(S)$}, the group $\pi_1 M_\phi$  has rank $2g+1$ and all minimal size generating sets are Nielsen equivalent. 

From above, when $A\in GL_2\BZ$ the group $G_A$ has rank $2$  exactly when there was some $v\in \BZ^d$ such that $\langle v, Av \rangle = \BZ^d$.  The Farey graph is the curve graph of $T^2$, and $\langle v, Av \rangle = \BZ^d$ exactly when $v,Av \in \mathcal F $ are adjacent, so in the Euclidean setting the analogue of  the rank part of our theorem in \cite{Biringerranks} still holds, and says that $\rank(G_A)=3$ if the translation distance of $A$ on $\mathcal F$  is at least two. The analogue of the Nielsen equivalence part is (a weaker version of) Question 1.

\subsection{Acknowledgements}

The author  is partially supported by NSF grant DMS 1611851. Thanks to Juan Souto for a helpful conversation, and to the referee for a number of useful comments.

%
%
%

\section{The proof}
\label {proof}
 We will first show that for a general $A \in GL_2\BZ$, there can be at most two Nielsen equivalence classes of $2$-element generating sets for $G_A$.   We'll then show that the conjugates of $\left(\begin {smallmatrix} 2 & 1 \\ 1 & 1 \end {smallmatrix}\right )$ are the only $A$ that realize this bound.  

\vspace{1mm}

The beginning of this argument overlaps with that of Levitt--Metaftsis \cite{Levittrank},  so we will just outline it and give citations when necessary. 
Suppose that $G_A$ is $2$-generated.  By \cite[Proposition 4.1]{Levittrank}, every minimal size generating set for $G_A$ is Nielsen equivalent to a generating set of the form
$$x=\left ( v, 0 \right), \ y=\left ( 0, 1 \right),  \text { where } v\in \BZ^2.$$

\comment{
Here, one can arrange that the second coordinates are $0,1$ using the Euclidean algorithm.  The proof of Theorem \ref {torus} shows that then $\left<v,Av\right>=\BZ ^ 2$, so if $w=av+bAv,$ then $\{(v,0),(w,1)\}$ can be changed to $\{(v,0),(0,1)\}$ by first multiplying $(w,1)$ by $(v ,0) ^ {- a}$, then conjugating $(v,0)$ to $(Av,0)$, then multiplying $(v ,0) ^ {- a}(w,1)$ by $(Av,0)^b$, and finally conjugating $(Av,0)$ back to $(v,0)$.}

Set $\CS_A=\{v\in \BZ^ 2\ | \ \left<v,Av\right>=\BZ ^ 2\}$.  Again by \cite[Proposition 4.1]{Levittrank}, if $v,v' \in\CS_A$, then $\{(v, 0),(0,1)\}$ and $\{(v', 0),(0,1)\}$ are Nielsen equivalent if and only if $v,v'$ lie in the same $\left<A\right>\times \BZ/2\BZ$-orbit on $\CS_A$, where $\BZ /2\BZ $ acts by $v\mapsto -v$.   

We now reinterpret this in terms of the Farey graph $\mathcal F$.   Recall from the  introduction that the vertex set of $\mathcal F$ consists of primitive elements of $\BZ ^ 2 $ up to negation, so can be identified with $\BQ \cup \{\infty\} $  through the map $$\mathcal F \longrightarrow \BQ \cup \{\infty\} ,  \ \ \pm \begin {pmatrix}a\\b \end {pmatrix}\mapsto \frac ab.$$ Below, we will regard $\BQ \cup \{\infty\}$ as a subset of $\BR \cong \partial_\infty \BH^2$, where $\BH^2$ is considered in the upper half plane model, and we will identify edges of $\mathcal F$ with the corresponding geodesics in $\BH^2$. (See all figures below.)  This embedding of $\mathcal F$ has some convenient properties. All edges of $\mathcal F$ separate $\BH^2 \cup \partial_\infty \BH^2$, and also $\mathcal F$, into two connected components. Every component of $\BH^2 \cup \partial_\infty \BH^2 \setminus \mathcal F$ is an ideal hyperbolic triangle, which we will call a \emph{complementary triangle} below. Finally, the action of $A \in GL_2\BZ$  on $\mathcal F $ is the restriction of its action on $\BH^2 \cup \partial_\infty \BH^2$ as a fractional linear transformation.  

Returning to the proof, vertices $v,w \in \mathcal F$ are adjacent if $ \langle v,w\rangle = \BZ^2$, so $\mathcal S_A $ is exactly the set of vertices in $\mathcal F$  that $A$  translates a distance of $1$.  Also, in the Farey graph we have identified primitive pairs up to negation, so the action of $\left<A\right>\times \BZ/2\BZ$ on $\CS_A$ is just the $A$-action on the corresponding set of vertices of $\mathcal F$. 
Define a \emph {1-orbit} of $A \circlearrowright \mathcal F$ to be an orbit all of whose points are translated a distance of $1$ by $A$.  
Theorem \ref {torus2} then becomes the following lemma.

\begin {lemma}\label {geometrytorus}
The action of $A \circlearrowright \mathcal F$ has a single $1$-orbit unless $A$ is conjugate to $\pm\left(\begin {smallmatrix} 2 & 1 \\ 1 & 1 \end {smallmatrix}\right )$, in which case it has two.

\end {lemma}

Fix a matrix $A \in GL_2\BZ $ and let $\ell$ be a $1 $-orbit of $A$.  Adding in edges connecting each $v\in \ell$ to $Av$, we will regard $\ell$ as an oriented path in $C (T ^ 2) $.  At each of its vertices $v $, a path $\ell$ has a \emph {turning number}, whose absolute value is one more than the number of Farey graph edges that separate the two edges of $\ell$ incident to $v $.  The turning number at $v $ is positive if the turn is counterclockwise when $\ell $ is traversed positively, and negative when the turn is clockwise.  (Remember that we are viewing $\mathcal F$ as a subset of the upper half plane in $\BR^2$.) When $v =\infty $, the turning number is just $A (v ) -A ^ {- 1} (v)$.  For instance, in Figure \ref {iteration} all turning numbers on the red 1-orbit are $3$, on the blue 1-orbit they are $-3$. 

When $A $ is orientation preserving, all the turning numbers on a given 1-orbit coincide.  On the other hand, if $A $ is orientation reversing then the turning numbers on a 1-orbit all have the same absolute value and alternate sign.  As $GL_2\BZ$ acts edge transitively on $\mathcal F$, any $1 $-orbit of $A $ may be translated to pass through $\infty, 0 $, which conjugates $A$ so that it has the form \begin {equation}\label {standardform}A=\begin {pmatrix} 0 & \epsilon \\1& x \end {pmatrix}, \ \ x\in\BZ, \ \epsilon =\pm 1. \end {equation}
When $A$ is as above, the turning number at $0 $ is $-\epsilon x$.  Checking eigenvalues, two matrices $\left (\begin {smallmatrix} 0 & \epsilon_i \\1& x_i \end {smallmatrix}\right)$, where $i=1,2$, are conjugate in $PGL_2\BZ$ if and only if $\epsilon_1 =\epsilon_2 $ and $|x_1|=|x_2|$.  This implies that the turning numbers of all the 1-orbits of a matrix $A$ have the same absolute value.

\begin{figure}[b]
\centering
\includegraphics{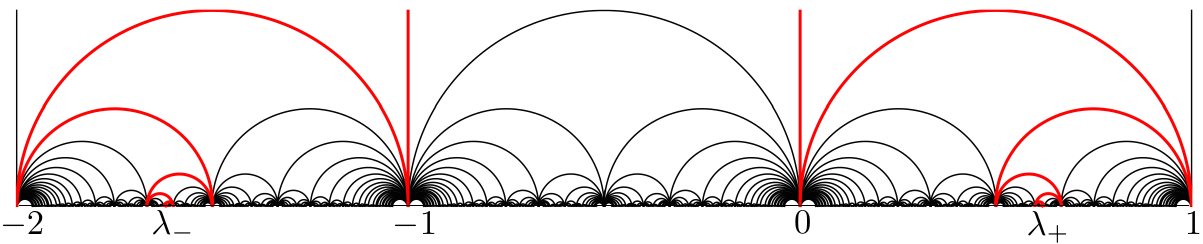}
\caption []{There is a single 1-orbit for the action $\left(\begin {smallmatrix} 0 & 1 \\ 1 & 1 \end {smallmatrix}\right )\circlearrowright \mathcal F$, on which the turning numbers alternate between $\pm 1$.
}
\label {reversingturning}
\end{figure}

It suffices to prove the lemma when $A= \left (\begin {smallmatrix} 0 & \epsilon \\1& x \end {smallmatrix}\right)$ as above.  Here, the conjugacy classes of $\pm\left(\begin {smallmatrix} 2 & 1 \\ 1 & 1 \end {smallmatrix}\right )$ correspond to the cases $\epsilon =-1,\ x=\pm 3$, so the goal is to prove that there are two 1-orbits in those cases, and one otherwise.  
\begin {itemize}
\item If $x=0$, then $A^2=\pm 1 $ and one can check directly that the only 1-orbit of $A$ is the edge connecting $\infty,0$.  
\item If $\epsilon =-1$ and $x=\pm 1$, then $A$ is orientation preserving and $A^3=\pm 1$.  Each of its $1$-orbits has turning number either $1$ or $-1$, so bounds a complementary triangle in $\mathcal F$.  But then $A$ is a rotation around the barycenter of this triangle in $\BH ^ 2 $, so this 1-orbit is the only one.
\item If $\epsilon =-1$ and $x=\pm 2$, then $A $ is parabolic.  Its 1-orbit has turning number $\pm 2$, so consists of all vertices in the $\mathcal F$-link of the fixed point of $A$.
\end {itemize}

When $A$ is hyperbolic, its 1-orbits are simple, biinfinite paths in $\mathcal F$ that accumulate onto the attracting and repelling fixed points $\lambda_+(A),\lambda_-(A)$.

\begin {itemize}
\item If $\epsilon =1$ and $|x|\geq 1$, then $A$ is hyperbolic and orientation-reversing.  The turning numbers on a 1-orbit $\ell$ alternate sign, so there is an edge of $\ell $ that separates  $\lambda_+(A)$ from $\lambda_- (A)$ in the upper half plane.  Any other 1-orbit would then have to intersect $\ell$, which is impossible, so $A$ has a single 1-orbit.  See Figure \ref {reversingturning} for an illustration of the case $\epsilon = 1 ,\ x=1$.
\item If $\epsilon=-1$ and $x=\pm 3$, then $A$ is orientation preserving, hyperbolic and conjugate to $\pm \left(\begin {smallmatrix} 2 & 1 \\ 1 & 1 \end {smallmatrix}\right )$.  When $x=3$, the orbits of $-1$ and $0$ are distinct, since they have opposite turning numbers (see Figure~\ref {iteration}).   Since the edge from $-1$ to $0$ in $\mathcal F$ separates the attracting and repelling fixed points of $A$, any 1-orbit of $A$ must pass through either $-1$ or $0 $.  So, the orbits of $\infty$ and $-1$ are the only 1-orbits.  The argument when $x=-3$ is similar. 

\end {itemize}

\begin{figure}[h]
\centering
\includegraphics{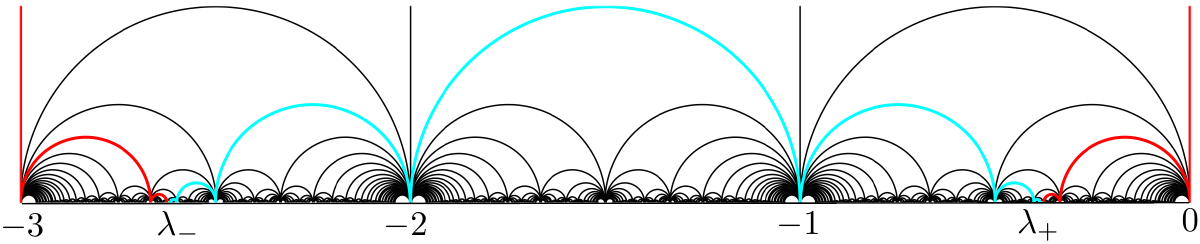}
\caption []{There are two orbits of $\left(\begin {smallmatrix} 0 & -1 \\ 1 & 3 \end {smallmatrix}\right )\circlearrowright \mathcal F$, its \emph {1-orbits}, on which the matrix acts as a translation by a distance of $1$.  The action is hyperbolic, with every forward orbit converging to $\lambda_+\approx-0.38$ and every backwards orbit converging to $\lambda_- \approx -2.62$. 
   Incidentally, the square of $\left (\begin {smallmatrix} 0 & 1 \\1& 1 \end {smallmatrix}\right)$ is a conjugate of $\left (\begin {smallmatrix} 0 & -1 \\1& 3 \end {smallmatrix}\right)$, which is why the vertex set of the 1-orbit in Figure \ref{reversingturning} is a translation of the union of the vertices of the two 1-orbits above. }
\label {iteration}
\end{figure}

It remains to deal with the case $\epsilon=-1$, $|x|\geq 4$, in which case $A$ is again orientation preserving and hyperbolic.  We claim that any biinfinite path $\ell$ whose turning numbers are all at least $3$ in absolute value is a geodesic in $\mathcal F$, and that if the turning numbers are all at least $4 $ in absolute value then $\ell$ is the unique geodesic in $C (T^2)$ connecting its endpoints.  This will imply that when $|x|\geq 4$, the matrix $A$ has only a single 1-orbit.

\begin{figure}[b]
\centering
\includegraphics{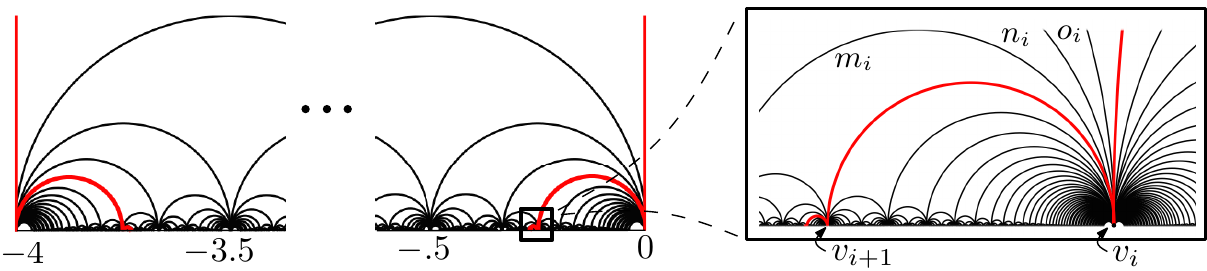}
\caption []{When $A=\left(\begin {smallmatrix} 0 & -1 \\ 1 & 4 \end {smallmatrix}\right )$, there is a single 1-orbit for the action $A \circlearrowright \mathcal F$, which is the unique geodesic connecting the attracting and repelling fixed points of $A$ in $\partial_\infty \BH ^ 2 $.}
\label {third}
\end{figure}

So, suppose that $\ell=(v_i)$ is a biinfinite path in $\mathcal F$ whose turning numbers are all at least $3 $ in absolute value.  For each $i$, let $m_i$ be the edge of $\mathcal F$ incident to $v_i$ that lies between the edges $[v_{i-1},v_i]$ and $[v_i,v_{i+1}]$, and shares a complementary triangle of $\mathcal F$ with $[v_i,v_{i+1}]$, as in Figure \ref {third}.  Each $m_i$ separates $m_{i-1}$ from $m_{i+1}$, so by planarity all the $m_i$ are disjoint.  Two vertices $v_i$ and $v_j$, with $i<j$, are disjoint from and separated by all the edges $$m_{i+1}\, ,\, \ldots \, , \, m_{j-1}.$$
Any path from $v_i$ to $v_j$ must go through all of these edges, so must have length at least $|i-j|$.  Therefore, $\ell$ is a geodesic in $\mathcal F$.

Suppose now that all the turning numbers of $\ell=(v_i)$ are at least $4$ in absolute value.  Choose for each $i $ two more edges $n_i,o_i$ incident to $v_i$ that lie between $[v_{i-1},v_i]$ and $[v_i,v_{i+1}]$, as in Figure \ref{third}.  All the edges $m_i,n_i,o_i$ separate the forward and backward limits of $\ell$, so any geodesic $\gamma$ in $\mathcal F$ connecting these limits must pass through a vertex of each $m_i,n_i,o_i$.  As $\gamma$ cannot pass through all three of the non-$v_i$ vertices of $m_i,n_i,o_i$, it must pass through $v_i$, so $\gamma=\ell$.  Thus, $\ell$ is the unique geodesic in $\mathcal F$ connecting its endpoints. This concludes the proof of Lemma \ref {geometrytorus}, and thus the proof of Theorem \ref {torus2}.

\begin{remark}
The educated reader will note that some of the simple properties of $\mathcal F$ used above reflect (and probably inspired) deeper results about the curve complexes of higher genus surfaces. For instance, the argument used to prove that a path whose turning numbers are all at least $3$ in absolute value is a geodesic is a simple version of Masur-Minsky's bounded geodesic image theorem \cite{Masurgeometry2}.
\end{remark}


\textit{\textrm{
\bibliographystyle{amsplain}
\bibliography{bibrefs}
}}

\end{document}